\documentclass[a4paper, 12pt, makeidx10pt]{article}
\usepackage[all]{xy}
\usepackage[dvipdfmx]{graphicx,color}

\usepackage{amsmath}
\usepackage{amssymb}
\usepackage{latexsym}
\usepackage{amsfonts}
\usepackage{mathrsfs}
\usepackage{bbm}
\usepackage{enumerate}
\usepackage{fancybox}
\usepackage{bbold}

\newcommand{\End}{\mathrm{End}}

\newcommand{\Ind}{\mathrm{Ind}}

\newcommand{\Res}{\mathrm{Res}}

\newcommand{\Inf}{\mathrm{Inf}}

\newcommand{\biperm}{{\sf biperm}}
\newcommand{\bbk}{{\mathbbm{k}}}
\newcommand{\bbK}{{\mathbbm{K}}}

\newcommand{\MTwoMot}{{\underline{\sf Mack}_{\bbk}}}
\newcommand{\cohMTwoMot}{{\underline{\sf Mack}_{\bbk}^{\sf coh}}}

\newcommand{\smod}{\mathbf{mod}}

\newcommand{\DoubleSpan}{{\sf S\widehat{pan}}}

\newcommand{\brg}{B(G)}  
\newcommand{\brgk}{B_{\bbk}(G)}  
\newcommand{\brgq}{B_{\Q}(G)}  
\newcommand{\brgR}{B_{R}(G)}  
\newcommand{\ghrg}{\widetilde{B}(G)} 
\newcommand{\ghrgk}{\widetilde{B}_{\bbk}(G)} 
\newcommand{\ghrgq}{\widetilde{B}_{\Q}(G)} 
\newcommand{\cbrg}{B^{\mathrm c}(G)}  
\newcommand{\cbrgk}{B^{\mathrm c}_{\bbk}(G)}  
\newcommand{\cbrgK}{B^{\mathrm c}_{\bbK}(G)}  
\newcommand{\cbrgq}{B^{\mathrm c}_{\Q}(G)}  
\newcommand{\cbrgo}{B^{\mathrm c}_{\calO}(G)}  
\newcommand{\cghrgk}{\widetilde{B}^{\mathrm c}_{\bbk}(G)}  
\newcommand{\cghrgK}{\widetilde{B}^{\mathrm c}_{\bbK}(G)}  
\newcommand{\cghrgq}{\widetilde{B}^{\mathrm c}_{\Q}(G)}  

\newcommand{\cclg}{\mathrm{C}(G)}  

\newcommand{\Q}{{\mathbb{Q}}}

\newcommand{\Z}{{\mathbb{Z}}}
\newcommand{\calA}{\mathscr{A}}

\newcommand{\calO}{\mathscr{O}}

\newcommand{\mcalM}{\mathcal{M}}

\newcommand{\mcalP}{\mathcal{P}}

\newtheorem{theorem}{\bf Theorem}[section]
\newtheorem{corollary}[theorem]{\bf Corollary}
\newtheorem{proposition}[theorem]{\bf Proposition}
\newtheorem{lemma}[theorem]{\bf Lemma}
\newtheorem{example}[theorem]{\bf Example} 

\newtheorem{remark}[theorem]{\bf Remark} 

\newcommand{\proof}{(Proof)\ }
\newcommand{\qed}{$\blacksquare$}

\newcommand{\iden}{\mathrm{id}}  
\newcommand{\Image}{\mathrm{Im}}  
\newcommand{\inv}{\mathrm{inv}}  

\usepackage{diagbox}
\usepackage{here}
\usepackage{mathtools}
\usepackage[dvipdfmx]{graphicx,xcolor}
\usepackage{tikz} 
\usepackage{tikz-cd}

\usetikzlibrary{cd}

\title{\vspace{-3cm}\bf Crossed Burnside rings and cohomological Mackey $2$-motives}
\author{
ODA, Fumihito\footnote{This work was supported by JSPS KAKENHI Grant Number 19K03457.} \\
 odaf@math.kindai.ac.jp \\
Dept. Mathematics, Kindai University, 
\\
Higashi Osaka, Japan
\\[\medskipamount]
TAKEGAHARA, Yugen \\ yugen@mmm.muroran-it.ac.jp \\
Dept. Mathematics, Muroran Institute of Technology, \\
Muroran, Japan
\\[\medskipamount]
YOSHIDA, Tomoyuki \\ yoshidat@math.sci.hokudai.ac.jp \\
Hokkaido University, \\
Sapporo, Japan
}
\date{\today}

\begin{document}
\maketitle

\begin{abstract}
Balmer and Dell'Ambrogio introduced the pseudo-functor $\mcalP$ from the bicategory of $\bbk$-linear Mackey $2$-motives to the bicategory of $\bbk$-linear cohomological Mackey $2$-motives over a commutative ring $\bbk$. They showed that $\mcalP$ maps the general Mackey $2$-motives to the cohomological Mackey $2$-motives by using the ring homomorphism from the crossed Burnside ring of a finite group $G$ over $\bbk$ to the center $ Z \bbk G$ of group algebra $\bbk G$ (\cite[Theorem 5.3]{Balmer Dell'Ambrogio21+}). 
We study the behavior of motivic decomposition of cohomological Mackey $2$-motives as images by $\mcalP$ of motivic decomposition of Mackey $2$-motives. 
\end{abstract}

\setcounter{section}{0}

\section{Introduction}
Let $G$ be a finite group and let $\bbk$ be a commutative ring. 
Balmer-Dell'Ambrogio theory of Mackey $2$-functors and Mackey 2-motives can be applied to a very wide range of research subjects. 
The theory suggests that the decompositions of many of the categories listed below as abelian categories are all controlled by Mackey $2$-motives. 
Here are some examples:
\begin{itemize}
\item In representation theory, $\bbk G$-modules $\mcalM (G)=\bbk G\mbox{-}\smod$ for the group algebra $\bbk G$ of $G$ over $\bbk$.  
\item In representation theory, derived categories $\mcalM (G) = {\mathrm{D}}(\bbk G)$ of $\bbk G\mbox{-}\smod$.
\item In equivariant homotopy theory, homotopy categories of $G$-spectra $\mcalM (G) = {\mathrm{SH}}(G)$.
\item In noncommutative geometry, Kasparov categories $\mcalM (G) = {\mathrm{KK}} (G)$ of $G\mbox{-}C^{\ast}$-algebras.
\end{itemize}
In order to unify and decompose the various abelian categories $\mcalM (G)$ mentioned above, Balmer and Dell'Ambrogio came up with the notion of Mackey $2$-motives which was inspired by the idea of the plain $1$-category of pure motives in algebraic geometry originated by Grothendieck.
They figured out that what controls the decomposition of $\mcalM (G)$ is the motivic decomposition 
\[
  G
\simeq  
(G,e_1)\oplus (G,e_2)\oplus \,\cdots \, \oplus (G,e_n),
\]
where 
$e_i$ in an idempotent of the endomorphism ring $\End_{\bbk\DoubleSpan(G,G)}({\mathrm{Id}_G})$ of the identity $1$-cell ${\mathrm{Id}}_G : G \to G$ in the bicategory $\bbk\DoubleSpan (G,G)$ (\cite[Definition 7.1.7]{Balmer Dell'Ambrogio20}) of $\bbk$-linear Mackey $2$-motives, and proved that the ring $\End_{\bbk\DoubleSpan(G,G)}({\mathrm{Id}_G})$ is isomorphic to the crossed Burnside ring $\cbrgk$ (\cite[Theorem 7.4.5]{Balmer Dell'Ambrogio20}) introduced by the third author (\cite{Yoshida97}).
Moreover, they introduced a theory of cohomological Mackey $2$-functors and cohomological Mackey $2$-motives in \cite{Balmer Dell'Ambrogio21+}.
As one of the many wonderful results that will be a result from Balmer-Dell'Ambrogio theory in the future, they introduced pseudo-functor $\mcalP$ which sends the ordinary Mackey $2$-motive to the cohomological Mackey $2$-motive in \cite{Balmer Dell'Ambrogio21+} to analyze in detail the relationship between the ordinary Mackey $2$-motive and the cohomological Mackey $2$-motive.

The main purpose of this paper is to give further refinements (Theorems \ref{theorem:A refinement of BD's Theorem 5.3}, \ref{theorem:A refinement of BD's Theorem 5.3 char(K) is 0}, and \ref{theorem:A refinement of BD's Theorem 5.3 dvr}) of \cite[Theorem 5.3]{Balmer Dell'Ambrogio21+} for some $\bbk$.
The second purpose of this paper is to determine the primitive idempotents of the crossed Burnside ring $\cbrg$ of $G$ over $\Z$ (Theorem \ref{theorem:primitive idempotents of cbr}). Although it is not directly related to the above topics, we will describe the properties (Proposition \ref{proposition:zeta is a surjective ring hom.} and Corollary \ref{corollary:existence of comm. diagram}) of the the ring homomorphism from $\cbrgk$ to the center $Z\mu_{\bbk}(G)$ of the Mackey algebra of $G$ over $\bbk$ introduced by Bouc. To investigate all primitive idempotents of $\cbrgo$ over a complete discrete valuation ring $\calO$, we apply Green functor theory by Bouc to give a decomposition of $\cbrgo$ (Proposition \ref{proposition:decomposition of CBR Green functor bG.}).

The paper is organized as follows: Section \ref{Section:IP of CBR} is a recollection of definitions and basic results on the crossed Burnside ring. Section \ref{Section:IP of p-localCBR} is a recollection of Bouc's theory of Green functors. As an application of this, we give the decomposition of $\cbrgo$, essentially, it was discussed in \cite{Bouc03}.
Section \ref{Section:Motivic decomp.} describes the behavior of motivic decomposition of cohomological Mackey $2$-motives as images by the several pseudo-functors $\mcalP$ of motivic decomposition of Mackey $2$-motives. Especially, we provide a condition to determine which primitive idempotent of the crossed Burnside ring is included in kernel of $\rho_G$ introduced by Balmer and Dell'Ambrogio (\cite[Theorem 5.3]{Balmer Dell'Ambrogio21+}).

In this paper, we fix a finite group $G$ with the identity $\epsilon$ and $\bbk$ be a commutative ring with identity.

\section{Idempotents of a crossed Burnside ring}\label{Section:IP of CBR}

\subsection{Crossed Burnside rings} We recall a construction of a crossed Burnside ring of $G$ over $\bbk$ from \cite{Bouc03}, \cite{OdaYoshida01}, and \cite{Yoshida97}. 
Denote by $G^{\rm c}$ the set $G$, on which $G$ acts by conjugation. 
The {\it category of crossed $G$-sets} is the category of $G$-sets over $G^{\rm c}$ : a crossed $G$-set $(X, \alpha)$ is a pair consisting of a finite $G$-set $X$ (i.e. a finite set with a left $G$-action), together with a $G$-map $\alpha$ from $X$ to $G^{\rm c}$, and a morphism of crossed $G$-sets from $(X, \alpha)$ to $ (Y,\beta)$ is a $G$-map $f$ from $X$ to $Y$ such that $\beta\circ f=\alpha$. 
The {\it crossed Burnside group} $\cbrg$ is defined as the Grothendieck group of the category of crossed $G$-sets, for relation given by disjoint union decomposition. 
Denote by $[X, \alpha]$ the isomorphism class of the crossed $G$-set $(X, \alpha)$. If $(X, \alpha)$ and $(Y, \beta)$ are crossed $G$-set, then their product is the crossed $G$-set
$(X\times Y, \alpha. \beta)$, where $X\times Y$ is the direct product of $X$ and $Y$, with diagonal $G$-action, and $\alpha.\beta$ is the map from $X\times Y$ to $G^{\rm c}$ defined by $(\alpha.\beta)(x, y) =\alpha(x)\beta(y)$. 
This product on crossed $G$-sets clearly commutes with disjoint unions, hence it gives a product on the group $\cbrg$. 
This turns $\cbrg$ into a ring. 
We call it {\it crossed Burnside ring of} $G$. 
The identity element of this ring is $[\bullet, u_{\bullet}]$, where $\bullet$ is a $G$-set of cardinality one, and the map $u_{\bullet}$ sends the unique element of $\bullet$ to the identity of $G$. 
A transitive crossed $G$-set is isomorphic to $(G/H,m_a)$ for a subgroup $H\le G$ and a map $m_a$ from $G/H$ to $G^{\rm c}$ by $m_a(gH)={}^ga:=gag^{-1}$ for an element $a\in C_G(H)$. Let $\mcalP_G$ denote the set of pairs $(H, a)$ consisting of a subgroup $H$ of $G$ and an element $a\in C_G(H)$. 
The group $G$ acts by conjugation on $\mcalP_G$, and we denote by $[\mcalP_G]$ a set of representatives of $G$-orbits on $\mcalP_G$. 
If $(H, a) \in\mcalP_G$, we denote by $[H, a]_G$ or $[(G/H)_a]$ the isomorphism class of the crossed $G$-set $(G/H,m_a)$. 
It is well known that a set $\{[H,a]_G\mid (H,a)\in [\mcalP_G]\}$ forms a basis of $\cbrgk$ (\cite[(3.1.c)]{OdaYoshida01}, \cite[Corollary 2.2.3]{Bouc03}) over $\bbk$. The ring has Burnside ring $\brgk$ of $G$ over $\bbk$ with basis $\{[G/H]\mid H\in \cclg\}$, where $\cclg$ is a complete set of conjugacy classes of subgroups of $G$, as a subring (see Lemma (\ref{Lemmas: Commutative diagrams and embeddings})). We denote by $\brg$ (or $\cbrg$), if the base ring $\bbk=\Z$. 

\subsection{Some maps between $\brg$ and $\cbrg$}
We define a ring homomorphism $\alpha_G^{\bbk}:\cbrgk\to\brgk$ by
\[
[(G/U)_t] \to  [G/U]
\]
and define a ring homomorphism $\iota_G^{\bbk}:\brgk\to\cbrgk$ by
\[
[G/U]\mapsto [(G/U)_{\epsilon}].
\]
Since $\alpha_G^{\bbk}\circ \iota_G^{\bbk}=\iden_{\brgk}$,
the Burnside ring $\brg$ is identified with $\Image\iota_G^{\bbk}$.
Let 
\[
 \ghrgk:=\prod_{H\in \cclg}\bbk.
\]
There exists a ring monomorphism $\phi^{\bbk}:\brgk\to \ghrgk$ given by 
\[
 [G/U]\mapsto (\phi_H^{\bbk}(G/U))_{H\in \cclg},
\]
where $\phi_H^{\bbk}([G/U])=\inv_H((G/U))=\{gU\in G/U\mid H\le {}^gU\}$.
The ring homomorphism $\varepsilon_H^{\bbk}:\bbk C_G(H)\to \bbk$ with $H\le G$ given by
\[
 \sum \ell_ss\mapsto \sum \ell_s
\]
is called the {\it augmentation map} of group algebra $\bbk C_G(H)$ of the centralizer $C_G(H)$ of $H$ in $G$ over $\bbk$ (cf. \cite[Definition 3.2.9]{MiliesSehgal02}).

We define a ring homomorphism $\widetilde{\alpha}_G^{\bbk}:\cghrgk\to\ghrgk$ where the subring 
\[
 \cghrgk=\left(\prod_{H\le G} Z \bbk C_G(H)\right)^G
\]
consists of $G$-fixed elements of product ring $\prod_{H\le G} Z \bbk C_G(H)$ of the centers of $\bbk C_G(H)$, by 
\[
 (x_H)_{H\le G}\mapsto (\varepsilon_H^{\bbk}(x_H))_{H\in\cclg}
\]
and define a ring homomorphism $\widetilde{\iota}_G^{\bbk}:\ghrgk\to\cghrgk$ by
\[
 (y_H)_{H\in\cclg}\mapsto (\widetilde{y}_H)_{H\le G},
\]
where $\tilde{y}_H=y_K$ for a conjugate $K\in\cclg$ of $H$ in $G$.
Obviously, $\widetilde{\alpha}_G^{\bbk}\circ\widetilde{\iota}_G^{\bbk}=\iden_{\ghrgk}$.
For a subgroup $H\le G$ there is a ring homomorphism $\varphi_H^{\bbk}:\cbrgk\to \cghrgk$ defined by linearlity 
by 
\[
 \varphi_H^{\bbk}([D,s])=\sum_{gD\in (G/D)^H}{}^gs=\sum_{t\in G}\sharp\{gD\in (G/D)^H\mid {}^gs=t\}\cdot t.
\]
The {\it Burnside homomorphism} is defined by
\[
 \varphi^{\bbk}=(\varphi_H^{\bbk})_{H\in\cclg} :\cbrgk\to\cghrgk.
\]
For simplicity, we will often use symbols $\varphi,\,\phi,\,\alpha_G,\,\iota_G$ etc. for $\bbk=\Z$.
We provide the following two lemmas.

\begin{lemma}\label{Lemmas: Commutative diagrams and embeddings}
\begin{itemize}
 \item[{\rm (i)}] The diagrams
\begin{equation}\label{Eq:Commutative diagrams}
 \begin{tabular}{ccc}
 \xymatrix{
    \cbrgk \ar[r]^{\varphi^{\bbk}} \ar[d]_{\alpha_G^{\bbk}} 
       & \cghrgk \ar[d]^{\tilde{\alpha}_G^{\bbk}} \\
    \brgk \ar[r]_{\phi^{\bbk}} & \ghrgk
  }
&
&
 \xymatrix{
    \cbrgk \ar[r]^{\varphi^{\bbk}} 
       & \cghrgk \\
    \brgk \ar[u]^{\iota_G^{\bbk}}\ar[r]_{\phi^{\bbk}} 
       & \ghrgk \ar[u]^{\tilde{\iota}_G^{\bbk}}
  }
%
\end{tabular}
\end{equation}
are commutative.
\item[{\rm (ii)}] Let $x\in\cbrg$. If $\varphi^{\bbk}(x)=\tilde{\iota}_G^{\bbk}(y)$ for some $y\in\ghrgk$, then $\iota_G^{\bbk}\circ \alpha_G^{\bbk}(x)=x$.
\end{itemize}
\end{lemma}

\proof
The statement (i) is clear. 
We prove the statement (ii).
Since $\tilde{\alpha}_G^{\bbk}\circ\tilde{\iota}_G^{\bbk}=\iden_{\ghrgk}$, it follows from the statement (i) that
\[
 \varphi^{\bbk}\circ \iota_G^{\bbk}\circ\alpha_G^{\bbk}(x)=\tilde{\iota}_G^{\bbk}\circ\phi^{\bbk}\circ\alpha_G^{\bbk}(x)=\tilde{\iota}_G^{\bbk}\circ\tilde{\alpha}_G^{\bbk}\circ\varphi^{\bbk}(x)=\tilde{\iota}_G^{\bbk}\circ\tilde{\alpha}_G^{\bbk}\circ\tilde{\iota}_G^{\bbk}(y)=\tilde{\iota}_G^{\bbk}(y)=\varphi^{\bbk}(x). 
\]
This shows that $\iota_G^{\bbk}\circ \alpha_G^{\bbk}(x)=x$, completing the proof.
\qed

\subsection{Primitive idempotents}
The primitive idempotents of the crossed Burnside algebra follows from a theorem of Bouc (\cite{Bouc03}) or \cite{OdaYoshida01}. 
The primitive idempotents of $\brgq$ have been determined by Gluck (\cite{Gluck81}) and the third author (\cite{Yoshida83})  independently. They are indexed by the conjugacy classes of subgroups of $G$. We denote by $e^G_H\in \brgq$ the primitive idempotent indexed by $H\le G$.  
The primitive idempotents of the Burnside ring $\brg$ follows from a theorem of Dress (\cite{Dress69}, or \cite[Corollary 3.3.6]{Bouc00}). 
We denote by ${\rm C}^{\infty}(G)$ a complete set of the conjugacy classes of perfect subgroups of $G$.
We write $H=_GK$ to denote that a subgroup $H\le G$ is $G$-conjugate to $K$.
The set of elements
\[
 f^G_J=\sum_{H^{\infty}=_GJ,\,H\in\cclg}e^G_H,
\]
where $H^{\infty}$ is the smallest normal subgroup of $H$ for which the quotient is soluble, for $J\in {\rm C}^{\infty}(G)$, is the set of primitive idempotents of $\brg$ (\cite[Corollary 5.4.8 (Dress)]{Benson 1991}).

The ring homomorphism $\iota_G:\brg\to \cbrg$ provides a decomposition of the identity of $\cbrg$ as a sum of orthogonal idempotents $\iota_G(f^G_J)$, for $J\in {\rm C}^{\infty}(G)$. 
We will show that the idempotents $\iota_G(f^G_J)'s $ are all the primitive idempotents of $\cbrg$.

\medskip
The next theorem is one of the main results of this paper and is applied to the proof of Theorem \ref{theorem:A refinement of BD's Theorem 5.3}.

\begin{theorem}\label{theorem:primitive idempotents of cbr}
The set of elements $\iota_G(f^G_J)$, for $J\in\mathrm{C}^{(\infty)}(G)$, is the set of primitive idempotents of $\cbrg$.
\end{theorem}
\proof
Let $x$ be an idempotent of $\cbrg$. According to \cite[Corollary 7.2.4]{MiliesSehgal02}, $\Z C_G(H)$ with $H\le G$ contains only trivial idempotents, where $\phi(x)=\tilde{\iota}_G(y)$ for some $y\in \ghrg$. 
This, combined with Lemma \ref{Lemmas: Commutative diagrams and embeddings} (ii), shows that $\iota_G\circ \alpha_G(x)=x$. By this fact, we may identify $x$ with $\alpha_G(x)\in \brg$. Since the map $\alpha_G:\cbrg\to \brg$ is a ring homomorphism, it follows that $\alpha_G(x)$ is an idempotent of $\brg$. Consequently, the idempotents of $\cbrg$ are those of $\brg$. This completes the proof.
\qed

\medskip

We may denote by $\overline{f}^G_J$ a primitive idempotent of $\cbrg$ indexed by $J\in{\rm C}^{\infty}(G)$ from Theorem \ref{theorem:primitive idempotents of cbr} above.
We denote by 
\begin{equation}\label{Eq:Commutative diagrams over the rationals}
 \begin{tabular}{ccc}
 \xymatrix{
    \cbrgq \ar[r]^{\varphi^{\Q}} \ar[d]_{\alpha_G^{\Q}} 
       & \cghrgq \ar[d]^{\tilde{\alpha}_G^{\Q}} \\
    \brgq \ar[r]_{\phi^{\Q}} & \ghrgq
  }
&
&
 \xymatrix{
    \cbrgq \ar[r]^{\varphi^{\Q}} 
       & \cghrgq \\
    \brgq \ar[u]^{\iota_G^{\Q}}\ar[r]_{\phi^{\Q}} 
       & \ghrgq \ar[u]^{\tilde{\iota}_G^{\Q}}
  }
\end{tabular}
\end{equation}
the commutative diagrams (\ref{Eq:Commutative diagrams}) for $\bbk=\Q$.

Since $\varphi_1$ (\cite[(4.2)]{OdaYoshida01},\cite[2.3.1]{Bouc03}) is a ring homomorphism from $\cbrg$ to the center $Z\Z G$ of group ring $\Z G$, as a matter of course we see that $\varphi_1(1_{\cbrg})=1$. More precisely, we can obtain an element $x\ne 1_{\cbrg}$ of $\cbrg$, where $x$ gives the identity of $Z\Z G$ as the image of $\varphi_{1}^{\Q}$.
\begin{corollary}\label{Corollary : the image of a primitive idempotent of cbr by phi} \it
Let $\overline{f}^G_J$ be a primitive idempotent of $\cbrg$ with $J\in {\rm C}^{\infty}(G)$. Then
\[
 \varphi_1^{\Q}(\overline{f}^G_J)=\left\{
\begin{array}{ll}
1 & J=1, \\
0 & {\rm otherwise}.
\end{array}
\right.
\]
\end{corollary}
\proof
Let $H\le G$. Since 
\begin{align*}
 \varphi^{\Q}_1\iota^{\Q}_G(e^G_H) 
&=  \widetilde{\iota}_G^{\Q}\circ \phi^{\Q}_1(e^G_H)\\
&=\left\{ 
\begin{array}{ll}
1 & H=1, \\
0 & {\rm otherwise}
\end{array}
\right.\\
\end{align*}
from (\ref{Eq:Commutative diagrams over the rationals}),
we obtain 
\begin{align*}
\varphi^{\Q}_1(\overline{f}^G_J)
&=\varphi^{\Q}_1(\iota^{\Q}_G(f^G_J))\\
&=\sum_{H^{\infty}=_GJ,H\in\cclg}\varphi^{\Q}_1(\iota^{\Q}_G(e^G_H))\\
&=\left\{ \begin{array}{ll}  1 & J=1, \\ 0 & {\rm otherwise}.  \end{array} \right.
\end{align*}
\qed

Let $\bbK$ be a field of characteristic $0$. Suppose that $\bbK$ is big enough. The primitive idempotents of $\cbrgK$ have been determined by \cite{OdaYoshida01} and \cite{Bouc03}. They are indexed by the subgroup $H$ of $G$ and irreducible $\bbK$-charater $\theta$ of $C_G(H)$. Let $e_{H,\theta}$ be a primitive idempotent of $\cbrgK$. Then since $e_{H,\theta}$ is given by $e_{H,\theta}:=\varphi^{-1}(\widetilde{\epsilon}_{H,\theta})$ where $\widetilde{\epsilon}_{H,\theta}$ is a primitive idempotent of $\cghrgK$ (\cite[Theorem (5.5)]{OdaYoshida01}), we have following result which prepares for proving Theorem \ref{theorem:A refinement of BD's Theorem 5.3 char(K) is 0}.

\begin{lemma}\label{Lemma : the image of a primitive idempotent of cbr over K by phi} \it
Let $e_{H,\theta}$ be a primitive idempotent of $\cbrgK$. Then
\[
 \varphi_1^{\bbK}(e_{H,\theta})=\left\{
\begin{array}{ll}
e_{\theta} & H=1, \\
0 & {\rm otherwise,}
\end{array}
\right.
\]
where $e_{\theta}$ is a primitive idempotent of $Z\bbK C_G(H)$ {\rm (c.f. \cite[Theorem 2.22]{NagaoTsushima88})}.
\end{lemma}

\section{Idempotents of a $p$-local crossed Burnside ring}\label{Section:IP of p-localCBR}

Let $\calO$ be a complete discrete valuation ring of characteristic $0$, with residue field $k$ of characteristic $p > 0$, and field of fractions is big enough. Review some results from Bouc's theory to obtain a decomposition of the crossed Burnside ring $\cbrgo$ of $G$ over $\calO$, we summarize the basic properties of Green functors and various results on its decomposition. 

\subsection{Green functors}

A {\it Mackey functor} $M$ for $G$ with values in the category $\bbk\mbox{-}\smod$
 of $\bbk$-modules is a bivariant functor $M=(M_{\ast},M^{\ast})$ from the category of finite $G$-sets to $\bbk\mbox{-}\smod$, with the following two properties:
\begin{enumerate}
 \item 
Let $X$ and $Y$ be any finite $G$-sets, and let $i_X\,(\mbox{resp. } i_Y)$ denote the canonical injection from $X\,(\mbox{resp. } Y)$ into $X\sqcup Y$. Then the morphisms
\begin{align*}
  (M_{\ast}(i_X),M_{\ast}(i_Y))
& :M(X)\oplus M(Y)\to M(X\sqcup Y),\\
{\footnotesize
\left(
\begin{array}{r}
M^{\ast}(i_X)  \\
M^{\ast}(i_Y) 
\end{array}
\right)
}
& :M(X\sqcup Y)\to M(X)\oplus M(Y)
\end{align*}
are mutually inverse isomorphisms.
\item Let
$$
\xymatrix{
    X \ar[r]^{a} \ar[d]_{b} 
       & Y \ar[d]^{c} \\
    Z \ar[r]_{d} & W
  }
$$
\normalsize
be a pull-back diagram of finite $G$-sets. 
Then 
\[
 M_{\ast}(b)\circ M^{\ast}(a) = M^{\ast}(d)\circ M_{\ast}(c).
\]
\end{enumerate}

A {\it Green functor} $A$ for $G$ over $\bbk$ is a Mackey functor for $G$ over $\bbk$, endowed for any $G$-sets $X$ and $Y$ with $\bbk$-bilinear maps $A(X)\times A(Y) \to A(X\times Y)$ with the following properties:
\begin{enumerate}
 \item If $f:X\to X'$ and $g:Y\to Y'$ are morphisms of $G$-sets, then the squares 

\footnotesize
 \begin{tabular}{cc}
\xymatrix{
A(X)\times A(Y)\ar[r]^{\times} \ar[d]_{A_{\ast}(f)\times A_{\ast}(g)}
    & A(X\times Y) \ar[d]^{A_{\ast}(f\times g)}\\
A(X')\times A(Y')\ar[r]_{\times}
    & A(X'\times Y')
  }
&
\xymatrix{
A(X)\times A(Y)\ar[r]^{\times} 
    & A(X\times Y) \\
A(X')\times A(Y')\ar[u]^{A^{\ast}(f)\times A^{\ast}(g)} \ar[r]_{\times}
    & A(X'\times Y')\ar[u]_{A^{\ast}(f\times g)}
  }
\end{tabular}
\normalsize

\noindent
are commutative.
\item If $X$, $Y$ and $Z$ are $G$-sets, then the square 
\footnotesize
\[
\xymatrix{
A(X)\times A(Y)\times A(Z)\ar[r]^{\iden_{A(X)}\times (\times)} \ar[d]_{(\times) \times \iden_{A(Z)}}
       & A(X)\times A(Y\times Z) \ar[d]^{\times}\\
A(X\times Y)\times A(Z)\ar[r]_{\times}
& A(X\times Y\times Z) 
  }
\]
\normalsize
is commutative, up to identifications $(X\times Y)\times Z\simeq X\times Y\times Z\simeq X\times (Y\times Z)$.
\item 
If $\bullet$ denotes the trivial $G$-set of cardinality $1$, there exists an element $\varepsilon_A\in A(\bullet)$, is called the {\it unit} of $A$, such that for any $G$-set $X$ and for any $a\in A(X)$
\[
 A_{\ast}(p_X)(a\times \varepsilon _A)=a=A_{\ast}(q_X)(\varepsilon_A\times a)
\]
denoting by $p_X$ (resp. $q_X$) the projection from $X\times \bullet$ (resp. $\bullet\times X$) to $X$.
\end{enumerate}
If $X$ and $Y$ are $G$-sets, if $a\in A(X)$ and $b\in A(Y)$, set
\[
 a\times^{op}b=A_{\ast}
(t)
(b\times a)\in A(X\times Y),
\]
where $t$ is a $G$-map defined by $Y\times X\to X\times Y; (y,x)\mapsto (x,y)$.
Define the {\it center} $Z(A)$ of $A$ by 
\[
 Z(A)(X)=\{a\in A(X)\mid \forall Y,\forall b\in A(Y), a\times b=a\times^{op}b\}
\]
for a $G$-set $X$ (\cite[12.1]{Bouc98}). The functor $Z(A)$ is a sub-Green functor of $A$.
If $e\in Z(A)(\bullet)$ is an idempotent, we denote by $e\times A$ the subfunctor of $A$ defined for a $G$-set $X$ by 
\[
 (e\times A)(X)=e\times A(X)\subset A(X).
\]
Then $e\times A$ is a sub-Green functor of $A$, with $e=e\times \varepsilon_A\in (e\times A)(\bullet)$ as unit. We denote by $W(H)$ the quotient group $N_G(H)/H$ for a subgroup $H\le G$.

Let $R$ be a ring in which every prime devisor of $|G|$ is invertible, except for $p$  which is not invertible. 
The primitive idempotents of the Burnside ring $\brgR$ follows from a theorem of Dress (\cite{Dress69}, or \cite[Corollary 3.3.6]{Bouc00}). 
We write $O^{p}(G)$ to denote that the smallest normal subgroup of $G$ for which the quotient is $p$-group.
A group $J$ is {\it $p$-perfect} if $O^p(J)=J$. 
We denote by ${\rm C}^{p}(G)$ a complete set of the conjugacy classes of $p$-perfect subgroups of $G$.
The set of elements 
\[
 f^G_J=\sum_{O^{p}(H)=_GJ,\,H\in\cclg}e^G_H,
\]
for $J\in {\rm C}^{p}(G)$, is the set of primitive idempotents of $\brgR$ (\cite[Corollary 5.4.8 (Dress)]{Benson 1991}).

For notations such as $\Ind_{N}^{G}$, $\Inf_{N/Q}^N$, etc. we follows that of Bouc's book \cite{Bouc98}.

\begin{theorem}{\rm \label{theorem:Decomposision of a Green functor}\cite[Proposition 12.1.11]{Bouc98}}
Let $R$ be a ring in which every prime devisor of $|G|$ is invertible, except for $p$ which is not invertible. Let $A$ be Green functor for $G$ over $R$. Then there are isomorphisms of Green functors 
\begin{equation}\label{Equation:decomposition of Green functor DVR}
 A\simeq \bigoplus_{J\in {\rm C}^{p}(G)}f^G_J\times A
\end{equation}
and
\begin{equation}\label{Equation:decomposition of Green functor quotient group DVR}
 f^G_J\times A\simeq \Ind_{N_G(J)}^G\Inf_{W(J)}^{N_G(H)}\left(f^{W(J)}_{1}\times(\Res^{G}_{N_G(J)}A)^J\right).
\end{equation}
\end{theorem}

\subsection{Dress construction}
We summarize the Dress construction of Green functors introduced by Bouc (\cite{Bouc03b}). See also \cite{OdaYoshida04}. 

Let $S$ be a finite $G$-set. If $M$ is a Mackey functor for $G$ over $\bbk$, then the Mackey functor $M_S$ is the bivariant functor defined on the finite $G$-set $Y$ by
\[
 M_{S}(Y)=M(Y\times S).
\]
If $f:Y\to Z$ is a map of $G$-sets, then
\[
 (M_S)_{\ast}(f)=M_{\ast}(f\times \iden_{S}),\quad 
 (M_S)^{\ast}(f)=M^{\ast}(f\times \iden_{S}).
\]
Then $M_{S}$ is a Mackey functor for $G$ and $M_S(\bullet)\cong M(S)$.

A {\it $G$-monoid} is a monoid endowed with a left $G$-action by monoid automorphisms. A {\it morphism of $G$-monoids} is a $G$-equivariant monoid homomorphism. A {\it crossed $G$-monoid} is a pair $(S,\varphi)$, where $S$ is a $G$-monoid, and $\varphi:S\to G^{\rm c}$ is a morphism of $G$-monoids.

\begin{proposition}{\rm \cite{Bouc03b}\label{proposition: associated Green functors}}
Let $(S,\varphi)$ be a crossed $G$-monoid. If $A$ is a Green functor for $G$ over $\bbk$, let $A_{S}$ denote the Mackey functor obtained by the Dress construction from the $G$-set $S$. If $X$ and $Y$ are finite $G$-sets, defined a product map $\times_S : A_{S}(X)\otimes_{\bbk}A_{S}(Y)\to A_S(X\times Y)$ by
\[
 \forall a\in A_{S}(X),\,\forall b\in A_S(Y),\,
a\otimes b\mapsto a\times_S b=A(\sigma)(a\times b),
\]
where $\sigma:X\times S\times Y\times S\to X\times Y\times S$ sending 
$(x,s,y,s')$ to $(x,\varphi(s)y,ss')$. Moreover, denote by $\varepsilon_{A_S}$ the element $A_{\ast}(f)(\varepsilon_A)$ of $A(S)\cong A_S(\bullet)$, where $f$ is the map sending the unique element of $\bullet$ to the identity element of $S$.
Then $A_{S}$ is a Green functor for $G$ over $\bbk$.
\end{proposition}


\subsection{Decomposition of a crossed Burnside ring over $p$-local ring}

Let $X$ be a $G$-set. We denote by $b(X)$ the Grothendieck group of the category of $G$-sets over $X$ : a $G$-set $(Y,\alpha)$ over $X$ is a pair consisting of a finite $G$-set $Y$, together with a $G$-map $\alpha:Y\to X$, and a morphism of $G$-sets over $X$ from $(Y,\alpha)$ to $(Z,\beta)$ is a $G$-map $f$ from $X$ to $Y$ such that $\beta\circ f=\alpha$ (\cite[Section 9]{Yoshida90}, \cite[2.4]{Bouc98}). 
Let $(Y,\phi)$ be a $G$-set over $X$. If $f:X\to X'$ is a morphism of $G$-sets, then we put $b_{\ast}((Y,\phi))=(Y,f\circ \phi)$. If $f:X'\to X$ is a morphism of $G$-sets, then we denote by $b^{\ast}((Y,\phi))$ the pull-back $(Y',\phi')$ of $(Y,\phi)$ along $f$, obtained by filling the cartesian square
\[
\xymatrix{
Y'\ar[r]^{a} \ar[d]_{\phi'}
       & Y\ar[d]^{\phi}\\
X'\ar[r]_{f}
       & X .
  }
\]
If $E=(U,\phi)$ (resp. $F=(V,\psi)$) is a $G$-set over $X$ (resp. over $Y$), we denote by $E\times F$ the $G$-set $(U\times V,\phi\times \psi)$ over $X\times Y$. Then the product $\times$ can be extended by linearity to a product from $b(X)\otimes_{\Z} b(Y)$ to $b(X\times Y)$.

\begin{proposition}{\rm \cite[Proposition 2.4.3]{Bouc98}\label{proposition: Burnside Green functors}}
With those notations above, $b=(b_*,b^{\ast})$ is a Green functor for $G$ over $\Z$.
\end{proposition}

We have a Green functor which gives a crossed Burnside ring $\cbrgk$.
Set $\bbk b(X)=\bbk\otimes_{\Z}b(X)$ for a $G$-set $X$.
The next proposition follows from Proposition \ref{proposition: associated Green functors}.

\begin{proposition}{\rm \cite[Theorem 5.1]{Bouc03b}\label{proposition: CBR Green functor}}
Let $\bbk b=(\bbk b_*,\bbk b^{\ast})$ be the Burnside Green functor for $G$ over $\bbk$ and $G^{\rm c}:=(G^c,\iden_{G^{\rm c}})$ be the crossed $G$-monoid.
Then $\bbk b_{G^{\rm c}}$ is a Green functor for $G$ over $\bbk$ and $\bbk b_{G^{\rm c}}(\bullet)\cong \cbrgk$.
\end{proposition}

We have a decomposition of a crossed Burnside ring $\cbrgo$ of $G $ over $\calO$. 



\begin{proposition}\label{proposition:decomposition of CBR Green functor bG.}
Let $R$ be a ring in which every prime devisor of $|G|$ is invertible, except for $p$  which is not invertible. 
Let $\{f_J^G\mid J\in {\rm C}^{p}(G)\}$ the set of primitive idempotents of $\brgR$.
Let $R b_{G^c}$ be the Green functor for $G$ over $R$. Then there is an isomorphism of Green functors 
\begin{equation}\label{Equation:decomposition of CBR Green functor DVR}
R b_{G^c}\simeq \bigoplus_{J\in {\rm C}^{p}(G)}f^G_J\times R b_{G^c}
\end{equation}
and
\begin{equation}\label{Equation:decomposition of CBR Green functor quotient group DVR}
 f^G_J\times R b_{G^c}\simeq \Ind_{N_G(J)}^G\Inf_{W(J)}^{N_G(H)}\left(f^{W(J)}_{1}\times(\Res^{G}_{N_G(J)}R b_{G^c})^J\right).
\end{equation}
In particular, there are ring isomorphisms
\begin{equation}\label{Equation:decomposition of CBR over DVR}
\cbrgo\cong \bigoplus_{J\in {\rm C}^{p}(G)}\iota_G^{\calO}(f^{G}_J) \cbrgo \end{equation}
and 
\begin{equation}\label{Equation:decomposition of CBR of qotient group over DVR}
 \iota_G^{\calO}(f^G_J) \cbrgo\cong \iota_{W(J)}^{\calO}(f^{W(J)}_1) B^{\rm c}_{\calO}(W(J)). 
\end{equation}
\end{proposition}

\proof
By Theorem \ref{theorem:Decomposision of a Green functor} and Proposition \ref{proposition: CBR Green functor}, we have isomorphisms 
(\ref{Equation:decomposition of CBR Green functor DVR}) and (\ref{Equation:decomposition of CBR Green functor quotient group DVR}) of Green functors.
By evaluation at the trivial $G$-set of those isomorphisms of Green functors, we have isomorphisms (\ref{Equation:decomposition of CBR over DVR}) and (\ref{Equation:decomposition of CBR of qotient group over DVR}) of rings as with Corollary 5.7.6 of \cite{Bouc00}. 
\qed

A Bouc's $\calO$-algebra 
\[
 \calA(G):=\cbrgo\iota^{\calO}_G(f^G_1)
\]
which contains all primitive idempotents $e \in\cbrgo$ such that $\iota^{\calO}_G(f^G_1) e=e$ is introduced by Bouc (\cite[3.2.3]{Bouc03}). 
He has been determined the primitive idempotents of $\calA(G)$ (\cite[3.2.11]{Bouc03}). They are indexed by the $p$-blocks of $ZkG$. 
We denote by $i_G$ the primitive idempotent of ${\calA}(G)$ corresponding to a $p$-block $i\in ZkG$ for a group $G$. 
We have a decomposition
\[
 \cbrgo\cong \bigoplus_{J\in {\rm C}^{p}(G)}\calA(W(J))
\]
into the direct sum of ideals by (\ref{Equation:decomposition of CBR over DVR}) of Proposition \ref{proposition:decomposition of CBR Green functor bG.}.
The set of elements 
$i_{W(J)}\in \calA(W(J))$, for $J\in {\rm C}^{p}(G)$ and a $p$-block $i\in Z kW(J)$, is the set of primitive idempotents of $\cbrgo$. Moreover, we have an expression 
\begin{equation}\label{Eq : the decomp. of 1 of cbrgo to primitive idemp.}
 1=\sum_{J\in {\rm C}^{p}(G),\,i\in Z kW(J)\,:\,\mbox{\footnotesize $p$-block}}
i_{W(J)}
\end{equation}
as the sum of primitive idempotents of the crossed Burnside ring $\cbrgo$ of $G$ over a discrete valuation ring $\calO$.
By the construction of $\calA(W(J))$ for $J\in {\rm C}^{p}(G)$, we have the following lemma.

\begin{lemma}\label{Lemma : the image of a primitive idempotent of cbr over dvr phi} \it
Let $i_{W(J)}$ be a primitive idempotent of $\cbrgo$, for $J\in {\rm C}^{p}(G)$ and $p$-block $i\in Z k W(J)$. Then
\[
 \varphi_1^{\calO}(i_{W(J)})=\left\{
\begin{array}{ll}
i & J=1, \\
0 & {\rm otherwise.}
\end{array}
\right.
\]
\end{lemma}

\section{Images of a motivic decomposition by pseudo-functor $\mcalP_{\bbk}$}\label{Section:Motivic decomp.}
Recall the {\it pseudo-functor} $\mcalP_{\bbk}$ from the $\bbk$-linear bicategory $\MTwoMot :=(\bbk\DoubleSpan^{\sf rf})^{\flat}$ of {\it Mackey $2$-motives} (see \cite[Recollection 2.2]{Balmer Dell'Ambrogio21+}) to the $\bbk$-linear bicategory $\cohMTwoMot :=(\biperm^{\sf rf}_{\bbk})^{\flat}$ of {\it cohomological Mackey $2$-motives} (see \cite[Definition 4.18]{Balmer Dell'Ambrogio21+}).
Balmer and Dell'Ambrogio showed the following theorem.

\begin{theorem}{\rm \label{theorem:Mackey 2-motives to cohomological Mackey 2-motives}\cite[Theorem 5.3]{Balmer Dell'Ambrogio21+}}
For every finite group $G$, there is a well-defined surjective morphism of commutative rings $\rho_G:\cbrgk\to Z(\bbk G)$ sending a basis element $[H,a]_G$ to
$\sum_{x\in [G/H]}{}^xa$. 
The pseudo-functor $\mcalP_{\bbk}$ maps the general Mackey $2$-motive $\oplus_{i}(G_i,e_i)$ to the cohomological Mackey $2$-motive $\oplus_i(G_i,\rho_{G_i}(e_i))$, where $(G,0)\cong 0$ in both bicategories.
\end{theorem}

\begin{remark}\rm
The ring homomorphism $\rho_G$ above is same as $\varphi_1$ in \cite[(4.2)]{OdaYoshida01} and $z_1$ in \cite[2.3.1]{Bouc03}. 
The map $\rho_G$ is not only a surjective ring homomorphism, but also essentially a special case of the homonymous one studied in \cite[CH. 7.5]{Balmer Dell'Ambrogio20}. 
\end{remark}

Although not directly necessary for discussions in this paper, we would like to mention here the relation between the ring homomorphism $\rho_G$ and the Mackey algebra of $G$.

Th\'evenaz and Webb introduced the {\it Mackey algebra} $\mu_{\bbk}(G)$ of $G$ over $\bbk$ in \cite{Thevenaz Webb 95}. 
We denote by $Y_{\bbk}(G)$ the endomorphism algebra
\[
 \End_{\bbk G}(\bbk \Omega_G)
\]
of permutation $\bbk G$-module $\bbk \Omega_G$
generated by the left $G$-set $\Omega_G=\coprod _{H\le G}G/H$ (\cite[Proposition 12.3.2]{Bouc98}, \cite[Definition 2.9]{Rognerud15}).
Th\'evenaz and Webb also introduced the {\it cohomological Mackey algebra} $co\mu_{\bbk}(G)$ of $G$ over $\bbk$. They pointed out that $co\mu_{\bbk}(G)$ is isomorphic to $Y_{\bbk}(G)$ and there exists a natural projection 
\[
 p_{\bbk}:\mu_{\bbk}(G) \to Y_{\bbk}(G)
\] (\cite[Section 16]{Thevenaz Webb 95}). 
The ring homomorphism 
\[
 \zeta_{\bbk}:\cbrgk\longrightarrow Z \mu_{\bbk}(G); [L,a]_G\longmapsto 
\sum_{U\le G}\sum_{w\in [L\setminus G/U]}t^U_{L^{w}\cap U}w^{-1}aw r^U_{L^{w}\cap U}
\]
is introduced by Bouc (\cite[Proposition 4.4.2]{Bouc03}). 
Moreover, he showed that the existence of an isomorphism 
\begin{equation}
\iota_{\bbk}: Z \bbk G\longrightarrow Z co\mu_{\bbk}(G),
\end{equation}
the details of the mapping are due to Rognerud (\cite[Proposition 3.4]{Rognerud15}) : the ring isomorphism $\iota_{\bbk}$ is defined by
\begin{align*}
Z \bbk G\ni  z 
&=\sum_{x\in G}\lambda_xx\longmapsto \\
\iota_{\bbk}(z)
&=\sum_{H\le G}
\sum_{g\in [H\setminus G/H]}
\left(\sum_{x\in H}\lambda_{gx}\right)
\overline{
t_{H\cap {}^gH}^H c_{g,H\cap H^g} r_{H\cap H^g}^H}
\in Z co\mu_{\bbk}(G).
\end{align*}

\begin{proposition}\label{proposition:zeta is a surjective ring hom.}
The homomorphism $\zeta_{\bbk}$ is surjective.
\end{proposition}

\proof
The natural map $p_{\bbk}:\mu_{\bbk}(G) \to Y_{\bbk}(G)\cong co\mu_{\bbk}(G)$
induces a ring homomorphism
\[
 p'_{\bbk}:Z \mu_{\bbk}(G) \to Z co\mu_{\bbk}(G).
\]
Since the surjectivity of $\varphi^1_{\bbk}$ and the commutative diagram
\begin{equation}\label{Diagram:zeta}
\xymatrix{
\cbrgk \ar[r]^{\zeta_{\bbk}} \ar[d]_{\varphi^1_{\bbk}}
       & Z \mu_{\bbk}(G) \ar[d]^{p'_{\bbk}} \\
 Z \bbk G \ar[r]_{\iota_{\bbk}}
     & Z co\mu_{\bbk}(G) 
  }
\end{equation}
shows that $p'_{\bbk}$ is surjective, we see that $\zeta_{\bbk}$ is surjective
from the commutativity of (\ref{Diagram:zeta}) again.
\qed

We recall the Burnside Green functor $\bbk b$ as used in Proposition \ref{proposition: CBR Green functor}.
Bouc introduced that $\bbk b(\Omega_G\times \Omega_G)$ has a ring structure which is isomorphic to $\mu_{\bbk}(G)$ (\cite[Proposition 4.5.1]{Bouc98}). See also \cite[Remark 2.5]{Rognerud15}. 

\begin{corollary}\label{corollary:existence of comm. diagram} \it
There exists a commutative diagram 
\begin{equation}\label{Diagram:beta and zeta }
\xymatrix{
\cbrgk \ar[r]^{\zeta_{\bbk}} \ar[d]_{\varphi^1_{\bbk}}
    & Z \mu_{\bbk}(G) \ar[d]_{p'_{\bbk}} \ar[r]^{\beta}
      &Z {\bbk}b(\Omega_G\times \Omega_G) \ar[d]^{p_L}\\
Z \bbk G \ar[r]_{\iota_{\bbk}}
    & Z co\mu_{\bbk}(G) \ar[r]_{\phi}
      &Z Y_{\bbk}(G)
}\end{equation}
of algebras.
\end{corollary}
\proof
The left square is commutative by (\ref{Diagram:zeta}) and the right square is also commutative by (2) of \cite[Theorem 2.12]{Rognerud15}.
\qed

\medskip
We will return to the discussion of Mackey $2$-motive.
In case of $\bbk=\Z$, we have a refinement of Theorem \ref{theorem:Mackey 2-motives to cohomological Mackey 2-motives} by Balmer and Dell'Ambrogio.

 \begin{theorem}\label{theorem:A refinement of BD's Theorem 5.3}
For every finite group $G$, the pseudo-functor $\mcalP_{\Z}$ maps the Mackey $2$-motive 
 \begin{equation}\label{Eq : motivic decomposition of G in Mackey 2-motive}
 G\simeq (G,\overline{f}^G_{1})\oplus (G,\overline{f}^G_{J_2})\oplus \cdots \oplus (G,\overline{f}^G_{J_m})
\end{equation}
to the cohomological Mackey $2$-motive 
\begin{eqnarray}
 G
& \simeq & (G,\rho_G(\overline{f}^G_{1}))\oplus (G,\rho_G(\overline{f}^G_{J_2}))\cdots \oplus (G,\rho_G(\overline{f}^G_{J_m}))\\
& \simeq & (G,{1})\oplus (G,0)\oplus\cdots \oplus (G,0)\\
& \simeq & (G,{1}),
\end{eqnarray}
where $\{1,\,J_2,\,\ldots,\,J_m\}={\mathrm{C}}^{\infty}(G)$.
\end{theorem}

\proof
By an argument of motivic decomposition in the proof of \cite[7.5.4]{Balmer Dell'Ambrogio20} and Theorem \ref{theorem:primitive idempotents of cbr}, we have an equivalence (\ref{Eq : motivic decomposition of G in Mackey 2-motive}). 
Then we have the rest of equivalences of cohomological Mackey $2$-motive from \cite[Theorem 5.8]{Balmer Dell'Ambrogio21+} and Corollary \ref{Corollary : the image of a primitive idempotent of cbr by phi}. The proof of the theorem is complete.
\qed

\begin{example}[Alternating group $A_5$]\rm
Let $G$ be the alternating group $A_5$ of $5$-letters.
Since $\mathrm{C}^{\infty}(G)=\{1,G\}$, Theorem \ref{theorem:primitive idempotents of cbr} shows that $\{\overline{f}^G_1,\overline{f}^G_G\}$ is the set of primitive idempotents of $\cbrg$. 
Moreover, by an example \cite[6.5 (F)]{OdaYoshida01}, 
we explicitly determine those elements as follows.
\small
\begin{eqnarray}
\overline{f}^G_1
&=& [A_4,\epsilon]+[D_{10},\epsilon]+[S_3,\epsilon]-[C_3,\epsilon]-2[C_2,\epsilon]+[1,\epsilon], \\
\overline{f}^G_G
&=& [A_5,\epsilon]-[A_4,\epsilon]-[D_{10},\epsilon]-[S_3,\epsilon]+[C_3,\epsilon]+2[C_2,\epsilon]-[1,\epsilon]. 
\end{eqnarray}
\normalsize
It is easy to see that $\rho_G(\overline{f}^G_1)=1$ and $\rho_G(\overline{f}^G_G)=0$. 
Therefore the pseudo-functor $\mcalP_{\Z}$ maps the Mackey $2$-motive 
\begin{equation}
 G\simeq (G,\overline{f}^G_{1})\oplus (G,\overline{f}^G_{G})
\end{equation}
to the cohomological Mackey $2$-motive 
\begin{eqnarray}
 G
& \simeq & (G,\rho_G(\overline{f}^G_{1}))\oplus (G,\rho_G(\overline{f}^G_{G}))\\
& \simeq & (G,{1})\oplus (G,0)\\
& \simeq & (G,{1}).
\end{eqnarray}
\end{example}

Let $\bbK$ be a field of characteristic $0$ and suppose that $\bbK$ is big enough. 
We denote by ${\rm Irr}(\bbK G)$ the set of all irreducible characters of $G$.
Lemma \ref{Lemma : the image of a primitive idempotent of cbr over K by phi} gives a refinement of \cite[Theorem 5.3]{Balmer Dell'Ambrogio21+} for $\bbk=\bbK$.

 \begin{theorem}\label{theorem:A refinement of BD's Theorem 5.3 char(K) is 0}
For every finite group $G$, the pseudo-functor $\mcalP_{\bbK}$ maps the Mackey $2$-motive 
 \begin{equation}\label{Eq : motivic decomposition of G in Mackey 2-motive char(K) is 0}
 G\simeq \bigoplus_{H\in \cclg,\,\theta\in{\rm Irr}(\bbK C_G(H))}(G,e_{H,\theta})
\end{equation}
to the cohomological Mackey $2$-motive 
\begin{eqnarray}
 G
& \simeq &\bigoplus_{H\in \cclg,\,\theta\in{\rm Irr}(\bbK C_G(H))}(G,\rho_G(e_{H,\theta}))\\
& \simeq &\bigoplus_{\theta\in{\rm Irr}(\bbK G)}(G,\rho_G(e_{1,\theta}))\\
& \simeq &\bigoplus_{\theta\in{\rm Irr}(\bbK G)}(G,e_{\theta}).
\end{eqnarray}
\end{theorem}

Let $\calO$ be a complete discrete valuation ring of characteristic $0$, with residue field $k$ of characteristic $p > 0$, and field of fractions is big enough. 
Lemma \ref{Lemma : the image of a primitive idempotent of cbr over dvr phi}
gives a refinement of \cite[Theorem 5.3]{Balmer Dell'Ambrogio21+} for $\bbk=\calO$.

 \begin{theorem}\label{theorem:A refinement of BD's Theorem 5.3 dvr}
For every finite group $G$, the pseudo-functor $\mcalP_{\calO}$ maps the Mackey $2$-motive 
 \begin{equation}\label{Eq : motivic decomposition of G in Mackey 2-motive dvr}
 G\simeq \bigoplus_{J\in {\rm C}^{p}(G),\, i\in Z k W(J)\,: \,\mbox{\rm \footnotesize  $p$-block}}(G,i_{W(J)})
\end{equation}
to the cohomological Mackey $2$-motive 
\begin{eqnarray}
G
& \simeq &
 \bigoplus_{J\in {\rm C}^{p}(G),\, i\in Z k W(J)\,: \,\mbox{\rm\footnotesize $p$-block}}(G,\rho_G(i_{W(J)}))\\
& \simeq & 
  \bigoplus_{i\in Z k W(1)\,:\, \mbox{\rm\footnotesize $p$-block}}(G,\rho_G(i_{W(1)})) \\
& \simeq & 
  \bigoplus_{i\in Z k G\,:\, \mbox{\rm\footnotesize $p$-block}}(G,i).
\end{eqnarray}
\end{theorem}
 
\renewcommand{\refname}
\end{document}